\documentclass[preprint,12pt]{elsarticle}
\usepackage{amssymb}

\journal{Applied Mathematics and Computation}

\begin{document}

\begin{frontmatter}

\title{\Large\bf Categorical Pairs and the Indicative Shift}
\author{Louis
H. Kauffman\\ Department of Mathematics, Statistics \\ and Computer Science (m/c
249)    \\ 851 South Morgan Street   \\ University of Illinois at Chicago\\
Chicago, Illinois 60607-7045\\ $<$kauffman@uic.edu$>$\\Tel: 1-312-996-3066\\Fax:1-312-996-1491.}

\begin{abstract}
The paper introduces the notion of a {\em categorical pair}, a pair of categories $(C,C')$ such that 
every morphism in $C$ is an object in $C'.$  Arrows in $C'$ can express relationships between the morphisms of $C.$ In particular we show that by using a model of the linguistic process of
naming, we can ensure that morphisms $F$ in $C$ can have an indirect self-reference of the form
$a \longrightarrow Fa$ where this arrow occurs in the category $C'.$ This result is shown to complement
and clarify known fixed point theorems in logic and categories, and is applied to G\"{o}del's Incompleteness Theorem, the Cantor Diagonal Process and the Lawvere Fixed Point Theorem.
\end{abstract}

\begin{keyword}
category; categorical pair; $2$-category; indicative shift;  self-reference;
indirect self reference.
\end{keyword}

\end{frontmatter}

\section{Introduction}

The purpose of this paper is to introduce a categorical pattern  that  complements
the Lawvere Fixed Point Theorem.  We produce a construction for indirect self-reference
that applies directly both to situations in ordinary language and to G\"{o}del's Theorem on the 
incompleteness of formal systems. Our construction can be summarized very succinctly and so we 
begin the  paper  with a self-contained account of the construction, and then devote the rest of the paper to discussion about how this {\em indicative shift} \cite{K1,K2,K3,K4,K5,K6,K7,K8,K9,K10,K11} can be seen in a number of different contexts. The indicative shift is defined in Section 2. The shift formalizes
an operation on names that can also be regarded as an {\em expansion} of a name in the sense that
if ``A" is the name of A then the expansion E``A" refers to A``A", the result of appending the contents of the name to the name. Thus if we regard the name as pointing to its contents as in 
$$``A" \longrightarrow A$$ then the expansion of the name points to the concatenation of the contents with the name $$E``A" \longrightarrow A``A".$$ Self-reference results when one expands the name of the expansion operator. $$E``E" \longrightarrow E``E".$$ The arrow of reference occurs at a different level than the concatenations of names and their contents. In this paper there will be two categories $C$ and $C',$ where the morphisms in $C'$ are arrows between morphisms in $C.$ We refer to $C'$ as the second category or the higher category in the pair $(C,C').$
The arrows of reference are in the higher category $C'.$
\bigbreak

In the case of G\"{o}del's Theorem, the pattern is similar. One has a method to assign natural numbers
(G\"{o}del numbers) to formulas. Letting $g \longrightarrow F(u)$ denote the assignment of the G\"{o}del number $g$ to the formula $F(u)$ with free variable $u$, we let $\sharp g$ denote the G\"{o}del number
that is assigned to $F(g)$, the formula obtained by evaluating $F(u)$ at its own G\"{o}del number.
Thus given $$g \longrightarrow F(u),$$ we have $$\sharp g \longrightarrow F(g).$$ Indirect self-reference is obtained by starting with a formula of the form $F(\sharp u).$ Then we have
$$g \longrightarrow F(\sharp u)$$ and 
$$\sharp g \longrightarrow F(\sharp g),$$ where $g$ is the G\"{o}del number of $F(\sharp u).$
The formula $F(\sharp g)$ refers to its own G\"{o}del number. This is the key indirect self-reference
behind the G\"{o}del theorem which then proceeds to take $F(u) = NB(u)$ where the formula $NB(u)$
asserts (within a formal system $\cal{F}$) that there is no proof of the statement with G\"{o}del number
$u.$ ($NB$ is short for ``nicht beweis".) Using this method of creating indirect self-reference we get
$$\sharp g \longrightarrow NB(\sharp g),$$
a statement that asserts its own unprovability in the formal system $\cal{F}.$ If the formal system is
consistent (and capable of handling these representations of arithmetic), then the statement $NB(\sharp g)$ is true, but unprovable in  $\cal{F}.$
\bigbreak

The paper is organized as follows. In Section 2 we construct a categorical context for the indicative shift by considering a pair of categories $(C,C')$ where every morphism in the first category is an object in the second category. Arrows in the second category can be interpreted as references between arrows in the first category. In this sense the second category of the categorical pair defined Section 2 takes the place of a meta-language in a logical context. We prove two basic results about self-reference in this section that we call the First and Second Self-Reference Theorems. The First Self-Reference Theorem gives conditions under which an indirect self-reference can occur. \bigbreak

\noindent The formality of the indicative shift is as follows. Suppose that $a$ and $F$ are morphisms in 
$C$ such that the composition $Fa$ is defined. 
Suppose that $$a \longrightarrow F$$ is an arrow in $C'$. Then it is either given or constructed
(First and Second Self-Reference Theorems) that there is an morphism $\sharp$ such that
$$\sharp a \longrightarrow Fa.$$ This is the indicative shift. It follows that if 
$$g \longrightarrow F\sharp$$ then
$$\sharp g \longrightarrow F\sharp g,$$ producing the desired indirect self-reference.
\bigbreak

The Second Self-Reference Theorem assumes that the pair of categories $(C,C')$ is a $2$-category, and that $a$ can be composed with itself. We show that the indicative shift with $\sharp a = aa$ follows naturally from the properties of composition in a $2$-category.   We end Section 2 with an application of its ideas to an example of Raymond Smullyan. Smullyan's example is a miniature version of G\"{o}del's incompleteness theorem. In Section 3 we show how the Self-Reference Theorems apply to 
G\"{o}del's Theorem when that theorem is seen as the production of a statement that asserts its own unprovability in a given formal system. In Section 4 we discuss the relationship of these ideas with the Lawvere Fixed Point Theorem, and we discuss how the Lawvere Theorem relates to G\"{o}del's Theorem via its analog with the Cantor diagonal process. The difference between our categorical approach and that of the Lawvere Theorem is that we formalize indirect self-reference.
\bigbreak

 In Section 5 we discuss how these category ideas and the indicative shift apply to ordinary language. More work needs to be done in relating these formalisms to ordinary language, for it is in ordinary language that the line between the categories (between level and meta-level)  is easily erased.
Names of names are still names in ordinary language, and in the language of categories, objects and morphisms can become interchangeable. From the mathematical side one can approximate the situation of language by using higher categories or even reflexive categories (where ideally there is a $1$-$1$ correspondence between objects and morphisms) rather than the categorical pairs  of Section 2. 
We give an example of a reflexive category (in that every object is a morphism) by taking the generating arrows and the objects to be the arcs of an oriented knot diagram. Section 6 is an Epilogue that reviews and discusses the ideas and results of the paper. The paper ends with Section 7, a return to self-reference and a discussion of the nature of self-reference in the use of the word I.
\bigbreak

\section{The Indicative Shift}
The indicative shift defined in this section formalizes
an operation on names that can also be regarded as an {\em expansion} of a name in the sense that
if ``A" is the name of A then the expansion E``A" refers to A``A", the result of appending the contents of the name to the name. Thus if we regard the name as pointing to its contents as in 
$$``A" \longrightarrow A$$ then $$E``A" \longrightarrow A``A".$$ Self-reference results when one expands the name of the expansion operator. $$E``E" \longrightarrow E``E".$$
\bigbreak

In those contexts where one thinks of expanding a name to its contents it is convenient to use the symbol
$E$ for the shift operator. In this section we shall adopt the symbol $\sharp$ for the shift. When we use
the symbol $\sharp$ we are thinking of the shift at the point where a name is given to a contents.
At that point, there is an initial pointing of the name to the contents before the name is directly associated 
with the contents. Then a shift occurs where the name is associated with the contents and the abstract
name is associated with the reference to these contents. These points of language are further discussed in Section 5.
\bigbreak

The reader should recall that a {\em category} \cite{Mac} consists in a collection of {\em objects} and a collection of
{\em morphisms}. To each morphism $f$ there is associated an ordered pair of objects $(A,B).$ We write $f:A \longrightarrow B$ to denote the morphism and call $A$ the {\em domain} of $f$ and
$B$ the {\em codomain} of $f.$ Given morphisms $f:A \longrightarrow B$ and $g:B \longrightarrow C,$
there is a morphism $g \circ f :A \longrightarrow C,$ called the {\em composition}
of $f$ and $g.$ Composition of morphisms is associative. Every object $A$
comes equipped with an identity morphism $1_{A}$ whose composition (with $A$ in the role of 
domain or codomain) with another morphism does not affect that morphism. This is the complete
definition of a category.
\bigbreak

We take as  given that in a category one can say whether two objects are equal and whether two 
morphisms are equal. We wish to model situations where equality is replaced by reference.
We speak this way for motivation and use the word reference as it is used in ordinary language where one may say that the name of a person refers to that person, or that the title of a paper refers to the 
text or to the contents of the paper.  We wish to model situations where one distinguishes between the morphisms of a given category and certain patterns of reference that are seen among these morphisms at a second level. An example that we shall consider later is the 
reference of a G\"{o}del number to its corresponding decoded text.  
\bigbreak

Let there be given a category $C$ and suppose that the set of morphisms of $C$ are seen as the objects
in another category $C'.$ We shall call the morphisms in $C'$ {\em reference arrows} for the morphisms
of $C,$ and we shall call the pair of categories $(C,C')$ a {\em categorical pair}. We make no further
restrictions on a categorical pair other than that there are two categories $C,C'$ with the morphisms in the first category forming the objects in the second category. A categorical pair is not constrained to be a
$2$-category (definition given below).
\bigbreak

There is a notion of {\em $2$-category} (and of higher categories) \cite{Le,Mac}. A 
$2$-category is a categorical pair with extra structure. In the notation of our reference arrows, the extra structure is as follows. One may have $$\alpha: a \longrightarrow b$$ and 
$$\beta: d \longrightarrow e,$$ arrows in $C'$ where
$ad$ and $be$ are both legal compositions in the base category $C$ of the categorical pair.
Then it is natural that there should be a referential arrow 
$$\alpha \circ_{0} \beta: ad \longrightarrow be,$$ usually called horizontal composition of these arrows in $C'.$ 
\bigbreak

Along with horizontal composition we have {\em vertical composition} which is simply the given composition of arrows in $C'.$ We can denote vertical composition by $\circ_{1}.$ Thus if
$$\alpha: a \longrightarrow b, \gamma: b \longrightarrow c$$
then $$\gamma \circ_{1} \alpha: a \longrightarrow c.$$
Now suppose that we have two possible vertical compositions 

$$\alpha: a \longrightarrow b, \gamma: b \longrightarrow c$$
$$\beta: d \longrightarrow e, \delta: e \longrightarrow f$$
where $ad$, $be$ and $cf$  are each legal compositions in the base category $C$ of the categorical pair.
Then it is natural to demand the compatibility
$$(\alpha \circ_{0} \beta) \circ_{1} (\gamma \circ_{0} \delta ) = (\gamma \circ_{1} \alpha) \circ_{0} (\delta  \circ_{1} \beta).$$
A categorical pair $(C,C')$ that satisfies this compatibility (called the {\em interchange law}) is called a 
{\em $2$-category}.  
\bigbreak

In the discussion below we shall formulate three self-reference theorems, one just using categorical pairs and the other two assuming that the categorical pair is a $2$-category. For most applications the reader will only need to assume given a categorical pair.
\bigbreak

Consider a categorical pair $(C,C').$  Let $a$ and $b$ be morphisms in $C$ and let
$a \longrightarrow b$ be a morphism in $C'$  with domain $a$ and codomain $b.$ Remember that
while $a$ and $b$ are morphisms in the initial category $C,$ they are objects in the referential category $C'.$ We call this arrow a {\em reference} from 
$a$ to $b.$ 
\bigbreak

We assume that for every object $X$ in $C$ there  is a morphism $\sharp_{X}: X \longrightarrow X.$ Thus each object in $C$ is indexed with a special morphism $\sharp_{X}$ to itself that need not be the identity morphism. Given an arbitrary morphism $a$ in $C$ we then have the compositions $a \sharp$ and $\sharp a,$
morphisms obtained by composing with $\sharp_{X}$ where $X$ is either the domain or codomain of $a.$ In this way, we will denote composition with $\sharp_{X}$ without explicitly writing $X.$
\bigbreak

We shall say that an arrow $a \longrightarrow b$ in the category $C'$ is a {\em composable reference} if the composition $ba$ is defined.
\bigbreak

We assume that the special morphisms $\sharp_{X}$ in $C$ have the the following 
property:
\bigbreak

\noindent {\bf The Indicative Shift.} {\em If $a \longrightarrow b$ is a composable reference arrow in 
$C',$  then there is an associated
reference arrow $\sharp a \longrightarrow ba.$ Here $\sharp a$ and $ba$ denote the compositions of these morphisms in the initial category $C$.} 
\smallbreak

\noindent A categorical pair with these properties is called a {\em referential pair.}
\bigbreak

\noindent {\bf Remark.} Note that since $C'$ is a category whose objects are the morphisms of $C,$ then given any morphism $a:A \longrightarrow B$ ($A$ and $B$ are objects in $C$), there is an identity morphism $a \longrightarrow a$ in
$C'$. Unless $A = B$ in $C$, this arrow will not be composable and so the shift will not apply to it.
When $A=B$, then we have a shift to $\sharp a \longrightarrow aa.$ Because we demand composability
of reference in order to have the indicative shift, not all morphisms in $C'$ can be shifted.  If the base 
category $C$ has only one object, or if all morphisms have the same domain and codomain, then any reference arrow can be shifted. In many examples, this simple circumstance is satisfied. Note also 
that we may have $ a \longrightarrow b$ composable with 
$$a: A \longrightarrow B$$ and 
$$b: B \longrightarrow Z$$
 so that the shift $\sharp a \longrightarrow ba$
exists, but no further shift is possible unless $Z = A.$
\bigbreak

\noindent {\bf First Self-Reference Theorem  (SRT1).} Let $(C,C')$ be a referential pair. Let $F$ be any morphism in the category
$C$, and assume that there is a composable reference $g \longrightarrow F \sharp$ in the category $C'.$  Then there exists a morphism $h$ in $C$ and a reference arrow in $C'$ such that 
$h \longrightarrow Fh.$
\smallbreak

\noindent {\bf Proof.} We are given  $$g \longrightarrow F \sharp$$  in $C'.$ This means that the codomain $X$ of $g$ is the domain of $F\sharp$, and this is the same as the domain of $F,$ and that 
$\sharp = \sharp_{X}.$  
Then apply the indicative shift and obtain:
$$\sharp g \longrightarrow F \sharp g.$$
Thus with  $h =  \sharp g,$ we have 
$h \longrightarrow Fh.$
This completes the proof. //
\bigbreak

\noindent {\bf Remark.} We interpret an arrow of the form $h \longrightarrow Fh$ as a model of an expression $Fh$ that is talking about (in the internal language of compositions in $C$) its own
``name" (which is the morphism $h$ from the point of view of the category $C'$).
\bigbreak

\noindent {\bf Remark.} Note that, in the above proof,  if $F:X \longrightarrow X$ (i.e. if $F$ was a morphism from an object to itself) then we could take the identity morphism  in $C'$ for
$F \sharp$, 
$$F \sharp \longrightarrow F \sharp$$ since then $F \sharp$ would be composable with itself.
Then the indicative shift produces
$$\sharp F \sharp \longrightarrow F \sharp F \sharp$$
and we have $a = \sharp F \sharp$ with $a \longrightarrow Fa.$ 
\bigbreak

\noindent {\bf Remark.}  Suppose that $g:X \longrightarrow X$ in $C$ and that $F:X \longrightarrow Y$
in $C.$ Then if we have $g \longrightarrow F$ in $C'$ then the indicative shift gives an infinite sequence of morphisms:
$$\sharp g \longrightarrow Fg$$
$$\sharp \sharp g \longrightarrow Fg \sharp g$$
$$\sharp \sharp \sharp g \longrightarrow Fg \sharp g \sharp \sharp g$$
and continuing in this fashion in the pattern
$$\sharp^{n} g \longrightarrow Fg \sharp g \sharp^{2} g \sharp^{3} g \cdots g \sharp^{n-1} g.$$
\bigbreak

On the other hand, if  $g:Z \longrightarrow X$ where the object $Z$ is distinct from X, then we can be given that $g \longrightarrow F$ and one shift to $\sharp g \longrightarrow Fg$ is possible. But the sequence of shifts stops here since the composition of $g$ with itself or with $\sharp g$ is not given to exist. This is one of the reasons for formulating this shift in categorical terms. The properties of 
the base category $C$ determine limits or lack of limits on the recursion of reference that is implicit in the 
indicative shift.  
\bigbreak

\noindent{\bf Remark.} For the next Theorem we will concentrate on references $a \longrightarrow F$
where the domain and codomain of $a$ are identical, so that there is no limit on the recursion of the shift. We call a morphism $a$ in $C$ a {\em self morphism} if it has the form $a: X \longrightarrow X.$
We can regard the referential arrows of the category $C'$ as generalizations
(categorifications) of the {\em equality} of morphisms in the base category $C.$ If the referential arrows are themselves taken to be equalities then the indicative shift would state that if $a = b$ as self morphisms in $C,$ then $\sharp a = ab.$ In other words, in this degenerate form, 
we would have  $\sharp a = aa$ for all self morphisms $a$ in $C.$
\bigbreak

\noindent The First Self-Reference Theorem would then correspond to the following calculation.
If $$g = F\sharp$$ then $$\sharp g = F\sharp g.$$ Hence $$gg = Fgg.$$
 The reader will recognise that this is exactly the form of the  proof of
the Church-Curry Fixed Point Theorem for Lambda Calculus \cite{B}. See the Epilogue (Section 6 of this paper) for more discussion of this point. The Indicative Shift generalizes the Church-Curry Fixed Point Theorem to a context that encompasses indirect self-reference.
\bigbreak

We formulate a second self-reference theorem that is close to the flavor of
the lambda calculus. We assume that the pair $(C,C')$ is a referential pair that is moreover a $2$-category in which $\sharp$ is defined in a particular way.
For a self morphism $a$ in $C$ we define
$$\sharp a = aa,$$ as the composition of $a$ with itself. We call the $2$-category a {\em lambda pair} 
if this condition is met for each self  morphism $a.$
\bigbreak

\noindent {\bf Second Self-Reference Theorem  (SRT2).} Let $(C,C')$ be a $2$-category that is a 
lambda pair as defined above with $\sharp a = aa$ for each self morphism in $C.$
Then, given a composable reference arrow  $a \longrightarrow F$ in $C',$ there is a corresponding
morphism $\sharp a \longrightarrow Fa.$ With respect to this indicative shift we obtain indirect 
self-reference from any composable reference  $a \longrightarrow F\sharp$ by taking the corresponding
shift to $\sharp a \longrightarrow F\sharp a.$ Note that this morphism is the same as
$aa \longrightarrow Faa.$ This final conclusion is a direct generalization of the Church-Curry Fixed Point Theorem.
\bigbreak

\noindent {\bf Proof.} Suppose we have a morphism $$a \longrightarrow F$$ in $C'.$
Let $a \longrightarrow a$ be the identity morphism for $a$ in $C'.$ Then we have the 
horizontal composition of these two morphisms:
$$aa \longrightarrow Fa.$$ Note that the composition $Fa$ exists since there is an arrow from $a$ to $F$ in $C'.$ Hence we have, as desired, the shift morphism
$$\sharp a  \longrightarrow Fa.$$ The rest of the Theorem follows in the same pattern as the 
proof of $SRT1.$ //
\bigbreak

\noindent{\bf The Smullyan Categorical Pair.} An exercise related to G\"{o}del's Theorem due to Raymond Smullyan \cite{S} can be naturally formulated in terms of categorical pairs. In this case we only use the structure of categorical pairs. We do not apply the indicative shift, but the Smullyan example contains its own indirect self-reference. The first category $C$ consists in (as morphisms) all words in the alphabet $\{ \sim, P, R, [\, ,\, ] \}.$ where a word is any ordered  string of these symbols. such words include the empty word which is the identity morphism in this category. The category $C$ has a single object. Composition in $C$ consists in concatenation of strings. The objects in the 
second category $C'$ consist in strings $X$ in the alphabet $\{ \sim, P, R, [\, ,\, ] \}.$  Thus every morphism in $C$ is an object in $C'.$   Other than the identity arrows, the following types of arrow in $C'$ are allowed, where $X$ is an arbitrary string in that alphabet.
\begin{enumerate}
\item $PX \longrightarrow P[X]$
\item $\sim PX \longrightarrow \sim P[X]$
\item $RX \longrightarrow P[XX]$
\item $\sim RX \longrightarrow \sim P[XX]$
\end{enumerate}
\bigbreak
Of course, once we allow these arrows in $C'$, we allow a host of possible compositions such as
the composition of $PX \longrightarrow P[X]$ and $P[X] \longrightarrow P[[X]]$ to form
$$PX \longrightarrow P[[X]].$$
The reader will note that by substituting $R$ for $X$ in item $3.$ we obtain the indirect self-reference
$$RR \longrightarrow P[RR].$$ 
By substituting $\sim R$ for $X$ in item $4.$ we obtain the indirect self-reference.
$$\sim R \sim R  \longrightarrow \sim P[\sim R \sim R].$$  \bigbreak

Smullyan has an amusing interpretation of this formalism. He tells the story of a 
machine that prints strings from the category $C$ (he does not use categorical terminology, but we will describe it that way). {\em Only the special itemized arrows (above) in $C'$ are interpreted as restrictions and descriptions of the machine's actions.} For codomains of arrows in $C',$ $P[X]$ means {\em printablity of $X$} and  $\sim P[X]$means {\em unprintablity of X.}  The category $C'$ contains the semantics for the categorical pair, but it also contains many expressions that are not interpreted semantically with regard to the machine's actions.
\bigbreak

\begin{enumerate}
\item If the machine can print the string $PX$ then it can print the string $X$. In other words
$$PX \longrightarrow P[X]$$ means that $X$ is individually printable if the string $PX$ is printable.
\item If the machine can print the string $\sim PX$ then the string $X$ is not printable (as an isolated 
string) by the machine. The string $\sim P[X]$ means that $X$ (alone) is not printable. 
$$\sim PX \longrightarrow \sim P[X].$$
\item The printing of the string $RX$ means that $XX$ is printable.
 $$RX \longrightarrow P[XX].$$
\item The printing of the string $\sim RX$ means that $XX$ is not printable.
$$\sim RX \longrightarrow \sim P[XX].$$
\end{enumerate}
Thus we can interpret the Smullyan Machine in terms of the category $C'$ by saying that certain 
special morphisms in $C'$ are interpreted as  statements about printability.  Each of the special string types (lets us call them {\em interpretable} strings)  $\{PX, \sim PX, RX, \sim RX \}$ might be printable by the machine, and if printed, they each tell what the machine can further print. It is given that {\em whenever
the machine prints one of these special strings then it tells the truth.} We deduce that
the machine cannot print the string $$\sim R \sim R \, ,$$ for this string asserts its own unprintability.
Thus, while the Smullyan Machine always tells the truth when it prints an interpretable string, there are interpretable strings that are true but unprintable!
This  Smullyan categorical pair is an intriguing miniature version of G\"{o}del's Incompleteness Theorem, with  printablity replacing provability. 
\bigbreak

\noindent{\bf A Simplest Example.} Let $C$ be a category with one object  $O$, the identity morphism $1_{O},$ and one other morphism $\sharp : O \longrightarrow O.$ We can take
the morphisms in $C$ to be the set of strings  $\{ \,\, , \sharp, \sharp \sharp, \cdots \}$ including the empty string identified as $1_{O}.$
In $C'$ we take as given  the morphism $$ \,\, \longrightarrow \,\,$$ from the empty string to the empty string. Then the shift (represented by $\sharp$) produces sequentially:
$$ \sharp \longrightarrow \,\, ,$$
$$ \sharp \sharp \longrightarrow \sharp ,$$
$$ \sharp \sharp \sharp \longrightarrow \sharp \sharp \sharp .$$
Self-reference appears at the third departure from the empty string. After that we have
$$\sharp^{n} \longrightarrow \sharp^{n(n-1)/2}.$$
\bigbreak

\noindent{\bf A Next Simplest Example.} Let $C$ be any category with one object  $O$ and morphisms $F$ other than the identity morphism $1_{O}$ and the morphism $\sharp: O \longrightarrow O$ representing the indicative shift.  

In $C'$ we have the identity morphism $$ X  \longrightarrow X$$  for each morphism $X$ in $C.$ Then the shift produces
$$ \sharp X \longrightarrow XX.$$
This pattern is analogous to the pattern of reference (by repetition) in the Smullyan Machine.
In particular we have the self-reference $$\sharp \sharp \longrightarrow \sharp \sharp$$ not necessarily 
the same as the identity morphism in $C'$ for $\sharp \sharp.$
\bigbreak

\noindent {\bf The Universal Building Machine.} We can interpret the expansion operator E described at the beginning of this section as a universal building machine. Then ``X" designates a blueprint for
the construction of X. (Of course here we indulge in a hierarchy of names. Really X is the name of 
an actuality and ``X" is the name of the blueprint for constructing this actuality.) Then we have
$$E``X" \longrightarrow X``X",$$ meaning that the universal builder E takes the blueprint ``X" and 
produces the actuality X appended to a copy of its blueprint. The  higher categorical morphism is a morphism between the composition of the building machine and the blueprint and the composition of the 
actuality and its blueprint. The universal building machine will build itself when supplied with its own blueprint. $$E``E" \longrightarrow E``E".$$
\bigbreak

\noindent {\bf Remark.} In the examples we have given so far, the category $C$ can be replaced with 
a monoid of strings under concatenation. In this $C$ there are many morphisms and only one object. Nevertheless, we have formulated the results in this section
to include categories $C$ with more than one object and the possibility of morphisms between distinct objects in $C$ Certainly in linguistic and other referential situations there are many examples where a given entity has a multiplicity of references to it.
This is modeled in a more general  category $C$ by a morphism $F$ such that there is a multiplicity of morphisms
$a$ such that the the composition $Fa$ is defined. Under these circumstances the indicative shift still holds and we may obtain a multiplicity of indirect self-references in the form 
$\sharp a \longrightarrow F\sharp a.$ An interesting source of abstract categories to consider for examples is found by starting with any directed graph $G$ and making a category $C$ whose objects are the nodes of $G$ and whose morphisms are generated by one identity morphism for each node, one
$\sharp$ morphism for each node and all the edges in the graph are interpreted as morphisms between their initial and final nodes. We allow multiple edges and loops in the graph. 
\bigbreak

\section{G\"{o}del's Theorem}
In order to discuss  G\"{o}del's Incompleteness Theorem from the point of view of the indicative shift, we
first start with the more general situation of a formal language $L$ that is susceptible to G\"{o}del
numbering. The basic notions of formal language and G\"{o}del numbering are explained in many books on logic. The interested reader can consult \cite{S,M,NN,Still}. We will assume that 
the formal language $L$ has the capacity to make statements about natural numbers involving a free variable $x$  such as ``The natural number  $x$ is greater than 2.". We will denote statements involving a single free variable $x$ in the form $S(x).$ Such a statement gives rise to infinitely many specialized 
statements that may be either true or false by substituting specific numbers for $x.$ Thus one could write 
$$S(x) = ``x >2"$$ and $$S(3) = ``3 > 2."$$
It is also understood that one can substitute the name or reference to a specific number for $x$ in a
statement $S(x).$ Thus, instead of the numeral $3,$  we could substitute for $x$ the statement `` the first odd prime number'' written in the language of $L$ (assuming that $L$ is rich enough to express this notion).
\bigbreak

Given a formal system, one can set up  G\"{o}del numbering, a method
that associates a unique natural number to each formula or sequence of formulas in $L.$ We write
$g \longrightarrow S(x)$ to denote the G\"{o}del number $g$ that is associated with a formula $S(x)$
with free variable $x.$ At this point the arrow is just a notation to indicate the association of the G\"{o}del number with its corresponding formula. We assume that there is a well-defined notion for substituting
a G\"{o}del number $g$ into a formula $S(x)$ to obtain a new formula $S(g).$ The new formula no longer has a free variable. What is substituted has to be a specific expression for the G\"{o}del
number in the language $L.$ Otherwise one would obtain a collection of formulas $S(h)$, one for each way to express the number $g.$  Once this choice has been made, then $S(g)$ has a specific G\"{o}del number. In particular,
we can start with $S(x)$, obtain its G\"{o}del number $g$ and then further obtain the  G\"{o}del number
$h$ of the result of substituting $g$ into $S(x).$ {\em We shall let $\sharp g$ denote a formula in $L$ that 
describes the process of computing the G\"{o}del number
of the result of substituting $g$ into $S(x).$}   Thus $\sharp g$ is a formula that stands for (the computation of) the G\"{o}del number $h.$ We shall write 
$$g \longrightarrow S(x)$$
and 
$$\sharp g \longrightarrow S(g).$$
Note that in the second equation, we use $\sharp g$ rather than $h$ on the left side of the arrow.
It is understood that $\sharp g$ stands for $h.$ The reader should note that while we use the arrow notation, no categories have yet been defined.
\bigbreak

We continue the story of these substitutions. We can assume that the formal system $L$
is rich enough to express in its own language the operation that takes the $g$ to 
 $\sharp g$ and that whenever one writes a formula of the form $S(x)$ one can also write
the formula $S(\sharp x).$ Here, as in the previous paragraph, $\sharp x$ stands for a formula in $L$ that describes the process of computing the G\"{o}del number
of the result of substituting the G\"{o}del number of $S(x)$ into the free variable in $S(x).$ 
Under these circumstances, we have 
$$g \longrightarrow S(\sharp x),$$
and 
$$\sharp g \longrightarrow S(\sharp g).$$
 We need to distinguish clearly between G\"{o}del numbers  and expressions
in the language $L$ that refer to the construction of such numbers. We assume that the 
G\"{o}del numbers are written in a standard numeral form like $3$ and not expressed indirectly as in
``the smallest odd prime number". 
The formula $S(\sharp g)$ refers to its own G\"{o}del number and hence achieves 
an indirect form of self-reference. In this formula, $g$ is a number written in the language $L$ and 
$\sharp g$ is a formula in $L$ that is applied to $g.$ Thus  the expression
$\sharp g$ refers to the G\"{o}del number of the formula $S(\sharp g).$
\bigbreak

This background is a short description of how indirect self-reference is accomplished in the context
of proving G\"{o}del's Incompleteness Theorem.  The rest of the well-known
proof of the Incompleteness Theorem uses this form of indirect self-reference applied to a statement
$S(x) =\,\, \sim B(x)$ that informally says ``The statement whose G\"{o}del number is $x$ has no
proof in the formal system $L$." The making of such a statement within the formal system $L$
requires that $L$ be sufficiently expressive so that it can internally encode the notion of a proof.
Once this is accomplished, one uses the construction of indirect self reference as shown below.
$$g \longrightarrow \, \sim B(\sharp x)$$
and 
$$\sharp g \longrightarrow \, \sim B(\sharp g).$$
The final statement $\sim B(\sharp g)$ asserts its own unprovability in $L.$ {\em  If $L$ is
consistent}, one concludes, by reasoning that occurs outside $L$,  that  $\sim B(\sharp g)$ is not
provable within $L.$ Thus $\sim B(\sharp g)$  is a statement that is true but unprovable by
$L.$ If the formal system $L$ is consistent, then it is incomplete. It has long been assumed that known
formal systems for elementary number theory are consistent. Under this assumption, such systems are incomplete. 
\bigbreak

We now indicate the categories $C$ and $C'$ that will place G\"{o}del's
Theorem in our context. A caution to the reader: These categories do not prove the Incompleteness Theorem. The proof still depends upon the careful construction of a formal system $L$ as described above. We obtain a description of how the indirect self-reference in the structure of the Incompleteness Theorem can be seen in a categorical framework.
\bigbreak

Let the base category $C$ have a single object, call it $O.$ Generating morphisms in $C$, other than the identity morphism and a special morphism $\sharp$,  are formulas in $L$ that have less than or equal to one free variable and natural numbers expressed outside the system $L.$ The natural numbers
outside the system are candidates for G\"{o}del numbers and will be composed according to those 
expressions inside $L$ for which they are code numbers.
A formula without a free variable may or may not define a number (integer) in $L.$ We shall call a formula {\em numerical} if it designates an integer. It is assumed that the language $L$ has a special category of formulas that designate numbers directly. These will be called {\em numerals}. For example one might use $|||$ as the numeral for $3$ in the formal system $L.$ When a 
G\"{o}del number is substituted into a formula the number is translated to the corresponding numeral
in $L.$ We will use usual decimal notation for numbers outside the formal system. The coding method could depend upon the decimal system (as in the example at the end of this section) or it could just depend upon number theoretic properties (such as 
the unique decomposition of a natural number into prime factors).
Let $G$ and $H$ be numerical formulas and let $n$ and $m$ denote numerals in $L$. Let $S(x)$ and $T(y)$ denote formulas with one free variable. Composition in $C$ will primarily  correspond to substituting one formula into the free variable of another formula. We define (non-identity) compositions in $C$  as follows:
\begin{enumerate}
\item It is given that composition is associative. 
\item$ S  \circ T$ is a formal composition with no specified relation if $S$ has no free variable and
$S$ and $T$ are formulas in $L.$
\item$ S(y) \circ T(x) = S(T(x))$ for $S$ a formula with a free variable and $T$ a formula in $L.$
Note that if $T$ has a free variable, then so does $S \circ T.$
\item $S(x) \circ n$ = $S(n)$ whenever $n$ is a numeral in $L.$
\item $S(x) \circ g$ = $S([g])$ whenever $g$ is any natural number. Here $[g]$ denotes the numeral in
$L$ that corresponds to $g.$
\item If $G$ and $H$ are formulas that represent numbers, but are not themselves numerals then 
$G \circ H $ is a formal composition with no specified relation.
\item If $n$ and $m$ are numerals in $L$, then 
$n\circ m$ is a formal composition with no specified relation.
\item If $g$ and $h$ are G\"{o}del numbers such that $g$ is the  G\"{o}del number of 
a formula with one free variable $S(x)$ and $h$ is the G\"{o}del number of a formula $T$,
then $g\circ h$ is the G\"{o}del number of $S([h]).$ Here we distinguish between G\"{o}del numbers
outside the system $L$ and numerical expressions inside that system.
\item If $g$ is the G\"{o}del number of a formula with one free variable, then we define
$\sharp \circ g = \sharp g$ as above. That is, if $g$ is the G\"{o}del number of $S(x)$, then
$\sharp g$ stands for the G\"{o}del number of $S([g]).$ Otherwise we take $\sharp \circ g$ formally with no specified relation.
\item If $n$ is a numeral in $L$ that stands for the G\"{o}del number $g$ of a formula with one free variable, then $\sharp n$ stands for $\sharp g$ expressed in $L$ as $[ \sharp g].$ It is understood that
$\sharp n$ is an expression in $L$ that refers to a specific numeral in $L.$
\item $\sharp \circ S(x)$ is taken formally with no specified relation.
\item $g \circ \sharp$ is taken formally with no specified relation for any natural number $g.$
\item $T \circ \sharp$ is taken formally with no specified relation for any formula $T$ in $L.$
\item $S(x) \circ \sharp$ is taken formally with no specified relation, but note that 
$(S(x) \circ \sharp) \circ g = S(x) \circ (\sharp \circ g) = S(x) \circ (\sharp g) = S(\sharp g)$
when $\sharp g$ is numerical. Similarly, $(S(y) \circ \sharp) \circ x = S(\sharp x).$
\end{enumerate}

This defines the category $C.$ We then define admissible arrows in $C'$ to be the identity arrows and
arrows of the form $$g \longrightarrow F$$ where $F$ is a formula in $L$ with at most one free variable,
and $g$ is the G\"{o}del number of $F.$  If $S$ is a formula in $L$ that has no free variable and $S$ represents the  G\"{o}del number $g$ of $F,$ then we also allow the morphism
$$S \longrightarrow F. $$ Compositions of morphisms in $C'$ are formal with no specified relations other than associativity. This completes the definition of the categorical pair corresponding to a given formal language $L$ with G\"{o}del numbering. In this way, the construction of indirect self-reference in G\"{o}del's Incompleteness Theorem can be regarded as an application of the First Self-Reference Theorem (SRT1).  
\bigbreak

\noindent {\bf Remark.} We can place G\"{o}del's Theorem in the context of the Second
Self-Reference Theorem {\bf SRT2.} Regard G\"{o}del numbers $g$ as morphisms in a category by 
defining $g \circ h$, as above,  to be the the result of substituting $h$ in the free variable of the decoding of $g$ (if there is such a free variable). Then we see that $\sharp g = g \circ g$ is a concise description of the $\sharp$ operator
as we have defined it above. With this definition of composition of G\"{o}del numbers we have a category
$C$ and can construct that category $C'$ of arrows from G\"{o}del numbers to texts in the formal system
just as we did in the above paragraphs. Now, if $g \longrightarrow F$ where $F$ is the decoding of $g$, then by definition $gg = \sharp g$
is the G\"{o}del number of $Fg$ where $Fg$ denotes the result of substituting $g$ into the free variable in $F.$ Thus we have $gg \longrightarrow Fg$ as the horizontal composition of $g \longrightarrow F,$
and the identity arrow $g \longrightarrow g$ and  $gg \longrightarrow Bgg$ as the horizontal composition of $g \longrightarrow B\sharp$
and the identity arrow $g \longrightarrow g.$This gives exactly the $2$-categorical structure of the 
Second Self-Reference Theorem.  The reader should note that in the category $C$ we have both Godel numbers and formulas as morphisms. In the category $C'$ we have identity morphisms that carry numbers to numbers and formulas to formulas, but otherwise arrows in $C'$ carry numbers or representatives of numbers to formulas. 
\bigbreak

In this way we see clearly that the categorification of the Church-Curry Fixed Point Theorem that is 
implicit in the Second Self-Reference Theorem applies to G\"{o}del's Theorem, showing how the 
indirect self-reference central to the G\"{o}del construction comes from  changing an equality to an arrow.
\bigbreak

\noindent {\bf Example.} In this example we give a small formal language that has  G\"{o}del numbering
and use it to illustrate our categorical constructions. Let ${\cal L}$ denote a language with the following
alphabet $${\cal A} = \{ (, ), \sim, P, x, |, \sharp \}.$$ The words in  ${\cal L}$ are all possible strings of
these symbols, and we will interpret them in a way that is similar to the Smullyan Machine described in
the previous section. Accordingly, we let $X$ denote any finite string of symbols in this alphabet.
We interpret $x$ as a variable. 
We interpret $|$ as the number $1$, $||$ as the number $2$ and generally a string $||| \cdots ||$ with 
cardinality $n$ vertical slashes as the {\em numeral} $n.$ The interpretation of $\sharp$ will be explained below.
We will introduce G\"{o}del numbering and show how to produce a string that refers to itself. 
If $X$ is a G\"{o}del number of a statement in ${\cal L},$  we will take that statement to be the 
{\em referent} of $X.$ We will interpret $P(X)$ to assert the {\em printability} of the referent of $X$ and 
$\sim P(X)$ to assert the unprintability of the referent of $X.$ We will construct a statement that asserts its
own unprintability. 
Then we will point out how the categories $C$ and $C'$ are constructed for the language 
${\cal L}.$
\bigbreak

The G\"{o}del numbers assigned to the individual members of the alphabet are as follows.
$$1 \longrightarrow ($$
$$2 \longrightarrow \, )$$
$$3 \longrightarrow \,  \sim $$
$$4 \longrightarrow P$$
$$5 \longrightarrow x$$
$$6 \longrightarrow \,  |$$
$$7 \longrightarrow \sharp$$
The G\"{o}del number assigned to a string of signs from the alphabet is the ordered list of these corresponding digits, interpreted as a natural number in the decimal system. Thus we have
$$34152  \longrightarrow \,\, \sim P(x)$$ so that $g = 34152$ is the G\"{o}del number of
$\sim P(x).$ If we wish to insert $g$ into the free variable in $\sim P(x), $ we must translate
$g$ into the language of numbers in this formal system. This means that we 
replace $34152$ by $|||\cdots |$ where there are {\em thirty four thousand one hundred fifty two} vertical slashes in this numeral. We can write the abbreviation
$\sim P(g),$  but in fact the actual expression is very large indeed.
Thus we have 
$$34152  \longrightarrow \,\, \sim P(x)$$ and 
$$341666  \cdots 62 \longrightarrow \,\, \sim P(||| \cdots |)$$
where there are $34152$ slashes in the right hand side, and there are the same number of $6$'s in the 
G\"{o}del number on the left hand side. Thus $341666  \cdots 62$ is the G\"{o}del number of the formula 
obtained from $\sim P(x)$ by substituting its own G\"{o}del number for the free variable $x.$ In our previous terminology we would write that, given $g = 34152,$ we have
$\sharp g = 341666  \cdots 62$ where there are $34152$ repeated $6$'s in the second number.
\bigbreak

We can now explain the interpretation of $\sharp$ in the formal system ${\cal L}.$ If a formula contains
$\sharp x$ for the variable $x$ or $\sharp n$ for some numeral $n,$  then $\sharp$ is
interpreted as that function that assigns to one  G\"{o}del number another G\"{o}del number by decoding the first number and placing its G\"{o}del numeral in the free variable of the decoded formula. If there is no such free variable, then no action is taken and there is no interpretation of $\sharp ||| \cdots ||.$
In a full formal system for number theory, $\sharp$ would be an abbreviation for an algorithm written in the language of that theory. Here we have only an external interpretation for the meaning of $\sharp.$
In the situation of the full formal system one has both the internal algorithm and the external interpretation, and one understands that these two ways of looking at the process describe the same function.
\bigbreak

Returning to our formula $$\sharp 34152  = 341666  \cdots 62$$ where there are $34152$ repeated 
$6$'s in the second number, we note that one can read the function $\sharp$ directly in the decimal formalism. If the decimal number $n$  contains the digit $5$ replace every occurrence of that digit by 
cardinality $n$ consecutive $6$'s. In this way $\sharp n$ is explictly described as a way to insert the number $n$ into itself by using two levels of coding the number (decimal and slash/numeral). When 
we speak of $\sharp n$ we can think about it as a direct function on decimal numbers, or we can understand it via its definition through coding and decoding in relation to the formal system  ${\cal L}.$
\bigbreak

We are now in  a position to see directly the composition of G\"{o}del numbers defined in this section
for the category $C.$ This category has one object and each string in the formal system ${\cal L}$ is
a morphism. Along with this, all natural numbers in the decimal system are also morphisms in $C$. Note that it is given that all natural numbers written in vertical slash numerals are morphisms in $C.$ We distinguish between slash numerals and decimal numerals.
Here are some examples of compositions of strings and numbers viewed as morphisms in the category
$C.$ 

\begin{enumerate}

\item Let $n$ and $m$ be any two decimal natural numbers. We define $n \circ m$ to be the decimal number $n$ obtained by substituting cardinality $m$ consecutive occurrences of the digit $6$ in
for every instance of the digit $5$ in $n.$ If $n$ does not contain the digit $5$ then no action is required for the composition of $n$ and $m.$
For example $4152  \circ 3 = 426661.$ Note that
$$ 4152 \longrightarrow \,\, P(x)$$ and 
$$426661 \longrightarrow \,\, P(|||).$$
Thus, if $n$ and $m$ are are G\"{o}del numbers of formulas $A$ and $B$,  and $A$ has a free variable $x,$ then $n \circ m$ is the G\"{o}del number of the result of substituting $B$ into the free variable in $A.$
If $n$ is the (decimal) G\"{o}del number of $A,$ then we have that $\sharp n = n \circ n.$

\item We have that $341752 \longrightarrow \,\,  \sim P(\sharp x)$ so that 
$\sharp 341752 = 341766 \cdots 62$ with  $341752$ consecutive $6$'s.
Thus $$341766 \cdots 62 \longrightarrow \,\,\sim P(\sharp ||| \cdots ||)$$
where there are $341752$ consecutive vertical slashes in the expression on the right.
This expression, $\sim P(\sharp ||| \cdots ||)$ refers to its own G\"{o}del number.
Hence, if the formal system ${\cal L}$ tells the truth, then the formula 
 $$\sim P(\sharp ||| \cdots ||)$$
is not printable, since it says that it
is unprintable. Since it is unprintable we find that it does tell the truth, and truth and printablity are 
distinct for the system ${\cal L}.$
\end{enumerate}

\bigbreak

\section{Lawvere's Fixed Point Theorem}
This section is a brief discussion of Lawvere's fixed point theorem. We discuss how the Lawvere Theorem arises as a generalization of the the Cantor diagonal argument, and we illustrate the Lawvere 
Theorem in the category of sets. For the reader interested in seeing the full formulation of this Theorem
in cartesian closed categories, we refer him to \cite{L1,L2}.
\bigbreak

Lawvere's Theorem \cite{L1,L2} is a direct generalization of Cantor's diagonal argument.
Recall Cantor's argument. We work in the category of sets. Let
$[A,B]$ denote the collection of set theoretic mappings from $A$ to $B.$
Let $Z = \{0,1\}$ and note that a subset $A$ of a set $X$ can be regarded as a mapping 
$A:X \longrightarrow Z$ where the elements of the subset are those $x \in X$ such that 
$Ax = 1.$ 
\bigbreak

\noindent {\bf Cantor.} Cantor gave a proof that {\em there is no surjective mapping from $X$ to $[X,Z].$}
His proof goes as follows. Let $F:X \longrightarrow [X,Z] $ be any mapping.
Define a subset $C$ of $X$ by the formula $$Cx = \sim F(x)x$$ where it is understood that 
$\sim 0 = 1$ and $\sim 1 = 0.$  $C$ cannot be of the form $F(a)$ for any 
$a \in X.$ For if $C = F(a),$ then $F(a)x = \sim F(x)x$ for all $x \in X.$ Hence
$ F(a)a = \sim F(a)a.$ This is a contradiction since the negation $\sim$ has no fixed points.
From this Cantor concludes that for $X$ infinite we have a higher infinity for $[X,Z]$ and so a 
hierarchy of infinities: $$X < [X,Z] < [[X,Z], Z] < \cdots.$$
\bigbreak

\noindent {\bf Lawvere.} Lawvere turns this scenario on its head by considering a more general case where $Z$ could be other than the set of two elements. 
\bigbreak

\noindent {\bf  Lawvere's Fixed Point Theorem for Sets.} Let $Z$ and $X$ be sets. Suppose that there exists a function $$F:X \longrightarrow [X,Z] $$ that is surjective. Let
$$\alpha: Z \longrightarrow Z$$ be any mapping from $Z$ to itself.  Then $\alpha$ has a fixed point.
\smallbreak

\noindent {\bf Proof.} Define $C: X \longrightarrow Z$ by the 
formula $$Cx = \alpha (F(x)x).$$ Then by surjectivity of $F,$ we have $C= F(a)$ for some $a$ and 
consequently $$F(a)a = \alpha(F(a)a).$$ Hence any mapping $\alpha: Z \longrightarrow Z$
must have a fixed point. //  \bigbreak

\noindent {\bf Remark.} Note that if we define a diagonal mapping $$\Delta : X \longrightarrow X \times X$$ by $$\Delta(x) = (x, x),$$ then $Cx = \alpha (F(x)x) = \alpha(eval((F \times I)(\Delta(x))))$ where $I$ denotes the identity map on $X$ and $eval(F(x), y) = F(x)y.$Thus
$$C = \alpha \circ eval \circ (F \times I) \circ \Delta.$$ In this way the map constructed in Lawvere's Theorem can be seen to work in any category with products and a terminal object. The terminal object serves to define the notion of a ``point". If $A$ is an object in the category and ${\bf 1}$ is the terminal object, then 
a {\em point in $A$} is a morphism $t:  {\bf 1} \longrightarrow A.$ Surjectivity of
$F:X \longrightarrow [X,Z] $ then means that for every $g: X \longrightarrow Z$ there is a
$t:  {\bf 1} \longrightarrow X$ such that $g = F \circ t.$ The diagonal map is crucial to the general construction. See \cite{L1} page 316 for a discussion of these points.
\bigbreak

Lawvere's Fixed Point Theorem can be used to place Cantor's orginal argument in different contexts.
For example, let $Z = \{0,1,J \}$ where $\sim 0 = 1, \sim 1 = 0, \sim J = J.$ In this example we can interpret $Z$ as the set of values in Lukasiewicz three-valued logic \cite{Luka}. Then generalized subsets of $X$ are described by maps into $Z.$ In such a generalized set $D: X \longrightarrow Z$, elements of $x \in X$ are either definite members of $D$ ($D(x) = 1$), definite non-members of 
$D$ ($D(x) = 0$), or indeterminate with respect to $D$ ($D(x) = J$).   If we have a mapping $F:X \longrightarrow [X,Z].$
then we can define a new mapping by the formula $Cx = \sim F(x)x$ and we find that if 
$C = F(z),$ then $F(z)z = \sim F(z)z,$ and we conclude that it must be the case that $F(z)z = J.$
We cannot conclude that $C$ is not of the form $F(z)$ for any $z \in X.$ Thus Cantor's argument 
about higher cardinalities does not generalize to a set theory based on the three-valued logic.
\bigbreak

\noindent {\bf Return to Self-Reference}
Now return to our First Self-Reference Theorem.  In this context, for the Russell Set, $Rx = \sim xx,$  we would generalize to
a reference arrow
$$R \longrightarrow \sim \sharp.$$
Applying the shift, we obtain
$$\sharp R \longrightarrow \sim \sharp R.$$
Instead of a contradiction, we obtain a referential arrow from the $\sharp R$ to its negation.
By changing equality to reference we have avoided the paradox. This is exactly how such paradox is resolved in computer languages where the referential step is often interpreted as a step in a recursive process. Of course we do not assert that this recursion solves the paradox in its original context.
\bigbreak

We end this section with a discussion of G\"{o}del's 
Incompletness Theorem in the Lawvere context and its relationship with our treatement of G\"{o}del
in the context of the First Self-Reference Theorem for categorical pairs.
\bigbreak

\noindent {\bf G\"{o}del Revisited.} Here is a how G\"{o}del's Theorem is related to the Lawvere Fixed Point Theorem. Let $\{ \phi (n,x) | n = 1,2,3,...\}$ denote a list of all syntactically valid formulas involving 
a single variable $x$ in the formal system $L$ (as described in Section 3). Suppose that $L$ is strong enough to be able (by proving or invalidating) to determine the truth or falsehood of each particular 
formula $\phi (n, m)$ for all natural numbers $n$ and $m.$   We define a new formula by 
$$Cx = \sim \phi(x,x).$$ Assuming that the list of all formulas and the ability of the formal system to determine their truth or falsity  is complete, we then have  $Cx = \phi (N,x)$ for some natural
number $N.$ Thus we have $$\phi(N,x) = \sim \phi(x,x)$$ for each natural number $x$ and hence 
$$\phi(N,N) = \sim \phi(N,N).$$ Since negation has no fixed point in the standard logic of $L$, we conclude that any list that we make of statements for the system will be of necessity incomplete with respect to the notion of truth within the system in terms of provability. Provability within,  and truth from  outside the system are distinct under the assumption  that the system $L$ is consistent.
 \bigbreak
 
 When we describe G\"{o}del's Theorem this way it is clear that it can be seen as an application of the Lawvere Fixed Point Theorem. We simply take $F(x)y = \phi(x,y)$ and the patterns match.
 Note that in this form of G\"{o}del's Theorem we did not encode directly a statement that asserts its own unprovability. This approach to G\"{o}del via a diagonal argument sidesteps the issue of self-reference. and instead shows the contradictory nature of completeness. This is the difference between the approach to G\"{o}del via the First and Second Self-Reference Theorems and the Lawvere Fixed Point Theorem.
Using the Self-Reference Theorems we construct an abstract framework for the G\"{o}del numbering and the indirect self-reference that is in back of the incompleteness phenomenon, and we show that this phenomenon is directly related to the higher categorical step of shifting from equality to 
 arrow.
 \bigbreak 

\section{Ordinary Language}
In this section we consider an interpretation for the First Self-Reference Theorem  in terms of ordinary language.
In this interpretation the morphisms of category $C$ are all  texts in ordinary discourse and all referents for these texts. Thus we regard perceptions and objects in the world as corresponding to  texts in a language that encompasses the written and spoken languages that are commonly used. In this way, if I meet another person, that other person would be regarded as a text whose name I come to learn in the course of meeting him or her. Then if I meet person P (a text) and learn his or her name N then at the beginning of that process there is indicated an arrow from N to P.
$$N \longrightarrow P$$
but shortly thereafter, when the naming process is more complete, the text that is P has become modified to contain its name in a prominent place and the name has been shifted to indicate that it is a name of 
{\em that person}. In actual practice this process is the one that includes our ability to recognise a person P as {\em that person with the name N}. We indicate this shift of reference by the indicative shift of
Section 2.
$$\sharp N \longrightarrow PN.$$
In terms of our perception, a text P that has undergone this shift is now known to have the name N.
The name N appears in our representational space along with the (text representing) the person P.
\bigbreak

Thus we see that the notion of categorical pair and indicative shift is a model of the referential shift inherent in the naming and referring of texts in ordinary language and in language in a very general
context. 
\bigbreak

The First Self-Reference Theorem  then becomes a model for how self-reference occurs in language.
For we see that the simplest instance of the Theorem is the act of naming the shift operation $\sharp.$
$$M \longrightarrow \sharp $$
Let $M$ denote the name of the shift operation $\sharp .$  Then $M$ is the name of the linguistic ability to combine a name with the text to which that name refers. And we see that once that name of the shift is
itself shifted, then a self-reference occurs.
$$\sharp M \longrightarrow \sharp M.$$
The completion of the naming process for the process of naming is self-referential. When we refer to ourselves in language we refer to our own ability to make and complete the act of naming.
\bigbreak

Note how the rest of the First Self-Reference Theorem works in this context. If we have a reference
$$G \longrightarrow F \sharp,$$ this is a reference to a text $F\sharp $ that talks about the naming process. Shifting this reference we obtain
$$\sharp G \longrightarrow F \sharp G,$$ a naming of a text that discusses its own name.
\bigbreak

We see that in the context of ordinary language a correct modeling must be flexible enough to allow even more hierarchies of reference and, at the same time to allow all these hierarchies to work at the same level since in language the name of a name is still a name. We see therefore that the splitting into two categories $C$ and $C'$ can lead to higher splittings (higher categories) and if these categories are all to be seen at a level, one may need to consider categories of infinite height where every object is a morphism and every morphism is an object. We call such categories {\em reflexive} and hope that they will be useful in an extension of this work to problems in mathematics,  linguistics and philosophy.
\bigbreak

\noindent To clarify these last remarks, consider a sequence of categories
$$C\, C' \, C'' \, \cdots C^{(n)} \,  C^{(n+1)} \,  \cdots $$ where the objects in $C^{(n+1)}$ are the morphisms
in  $C^{(n)}.$ We shall say that the category $C^{(n)}$ is of {\em type $n.$} There are a number of competing definitions for the notion of {\em $n$- category} (recall our specific definition of $2$-category in Section 2). All $n$- categories are of type $n.$ All the pair constructions in this paper apply in the transition between $ C^{(n)}$ and 
$C^{(n+1)}.$ A {\em reflexive category} is at level $C^{\infty}$ where any finite
descent from morphism to object will reveal only further morphisms. 
\bigbreak

It might seem that a reflexive category would be a huge undertaking, requiring some sort of limiting 
construction from a hierarchy of categories. That this is not so is illustrated in Figure 1. Here we
show, at the top of the Figure, a morphism between two morphisms. If one were to draw a diagram of morphisms such that every morphism of the diagram occurred between two morphisms, then the diagram could be interpreted as describing a reflexive category whose objects are the morphisms depicted in the diagram, and whose morphisms are generated under composition by the morphisms shown in the diagram
(and the unwritten identity morphisms, one for each pictured morphism). In this way, certain special
diagrams can represent reflexive categories.
\bigbreak

In Figure 1 we depict such a diagram. We show a diagram $T$ of a trefoil knot and take the oriented arcs of that diagram to be morphisms in a category that we shall call the {\em Trefoil Category.}  Knot diagrams have just the right properties, as described above, to generate reflexive categories. 
The generating morphisms are the arcs in the diagram and we take the objects of the category to also be this set of arcs. A morphism begins at one arc and ends at another arc. Every morphism in this category is a morphism of morphisms. Knot diagrams are of independent interest as they are formalizations of projections of curves in three space and can be used to faithfully study the topology of curve embeddings in three dimensional space. For this purpose one usually takes the knot diagrams up to an equivalence relation generated by the graphical moves shown in Figure 2 (the Reidemeister moves).
It is not our purpose here to dwell on the theory of knots, but in fact this association of a category to a knot diagram can, in principle, be used to obtain topological information about the knot. We will treat this
aspect in a separate paper.
\bigbreak

We also illustrate in Figure 1 a diagram $T'$ that is not quite a knot diagram that has the same formal characteristic of generating a reflexive category.
Each arc is seen to be an arrow originating on one of the arcs and terminating on another.
If the reader examines the Figure, it will be apparent that we have a category with objects
$\{A,B,C\}$ and each of these objects is a morphism with
\begin{enumerate}
\item $A: C \longrightarrow B,$
\item $B: A \longrightarrow C,$
\item $C: B \longrightarrow A.$
\end{enumerate}
Compositions of these morphisms are available, so this category has more morphisms than it has objects, but it is certainly reflexive in that all its objects are morphisms. Reflexive categories of this 
sort can be associated with knots and links. We shall study them in a separate paper.  A second example is shown in Figure 1 with the link diagram $L.$ Here the associated reflexive category has two objects $A$ and $B$ that are also generating morphisms for the category. We have
\begin{enumerate}
\item $A: B \longrightarrow B,$
\item $B: A \longrightarrow A.$
\end{enumerate}
The distinct morphism/objects $A$ and $B$ are ``linked" categorically in that each plays the role of 
a morphism for the other. It is clear that this notion of linking is close to the way we speak of linking in 
ordinary language where a linkage of plans, ideas or persons involves how each is a process for the other. One reason for 
bringing in this example of a reflexive category in a section on ordinary language is that we see that 
the Trefoil Category and the Link Category  (and indeed the diagrams as mathematical structures) arise from the language of sketching of three dimensional forms. But also, when we translate these diagrammatic forms into the corresponding reflexive categories we see that the categories themselves contain patterns of mathematical/topological languager.  These are topics to be pursued elsewhere.
 
 \begin{figure}[htb]
     \begin{center}
     \begin{tabular}{c}
     \includegraphics[width=5cm]{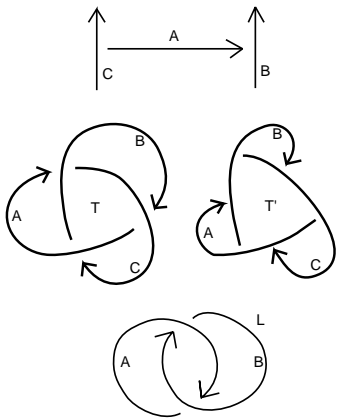}
     \end{tabular}
     \caption{\bf Trefoil Category and Link Category}
     \label{Figure 1}
\end{center}
\end{figure}

 \begin{figure}[htb]
     \begin{center}
     \begin{tabular}{c}
     \includegraphics[width=5cm]{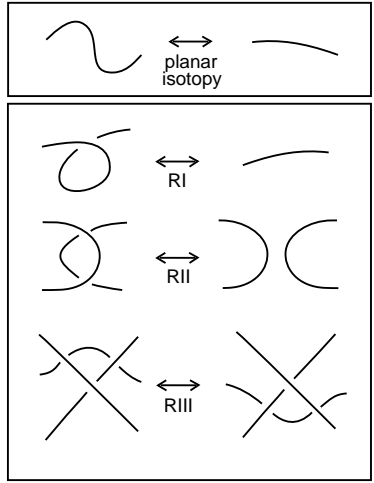}
     \end{tabular}
     \caption{\bf Reidemeister Moves}
     \label{Figure 2}
\end{center}
\end{figure}

 \section{Epilogue}
 Both the Self-Reference Theorems of this paper and the Lawvere Fixed Point Theorem
 come from generalizing the Cantor diagonal process, and both can also be seen as ways to generalize
 the Church-Curry Fixed Point Theorem.   In the Church-Curry Theorem
 we are given an algebra with a binary operation that is not  associative and an {\em axiom of reflexivity}  that states
 that {\em functions of a single variable expressed in that algebra can be named and regarded as elements of the algebra.} Thus in such an algebra $\Lambda$ one might define $G[x] =a((bx)x)$ as a function from the algebra to itself. One is then guaranteed that there exists an element $g$ such that for all $x$ in the algebra, $gx =  a((bx)x).$ This reflexive assumption of a correspondence between elements of the algebra and mappings of the algebra to itself is very strong. 
 \bigbreak
 
The simplest instance of this strength is the 
 Church-Curry Fixed Point Theorem which states that every element $F$ of $\Lambda$ has a fixed point in the sense that there is an $a$ in the algebra such that $Fa = a.$ The proof goes as follows.
Define $G[x] = F(xx)$ for all $x$ in $\Lambda.$ Then, by the axiom of reflexivity  there exists $g$ in $\Lambda$ such that $gx = F(xx)$ for all $x.$ Letting $x = g$ we obtain $gg = F(gg).$ So $gg$ is the fixed point for $F.$
\bigbreak

At the formal level, the Lawvere Fixed Point Theorem can be seen as a categorical generalization of the
$\Lambda$ algebra formalism $C[x] = \alpha(F(x)x)$ where it is known that such a $C$ must be represented algebraically by an element of the form $F(a)$ (the surjectivity hypothesis for $F$).
Then we have $F(a)x = \alpha(F(x)x)$ and consequently $F(a)a = \alpha(F(a)a),$ giving $\alpha$ 
a fixed point with a specific structure. The generality of the pattern allows it to be applied to many situations beyond the original Cantor argument. The application of the Fixed Point Theorem  to G\"{o}del's Theorem works best when we do not think of G\"{o}del's Theorem as depending on indirect self-reference.
\bigbreak

The First and Second  Self-Reference Theorems are generalizations of the Church Curry Fixed Point Theorem where we replace equality signs by  arrows of reference and we correspondingly generalize the operator
 $\sharp x = xx$ to an arrow of reference $$\sharp x \longrightarrow xx.$$
We then generalize the fundamental repetition operator $\sharp$ a notch further to the indicative shift
where, if
$$a \longrightarrow b$$
then 
$$\sharp a \longrightarrow ba$$
and the Church-Curry Fixed Point Theorem is transformed into  our First Self-Reference Theorem. In fact we could take the initial category $C$ to have one object and its morphisms the elements of the lambda algebra having either no free variable or a single free variable. Composition $ab$ of morphisms $a$ and $b$ is defined whenever $a$ has a free variable. Then $ab$ stands for the substitution of $b$ into the free variable in $a.$ With this we have both the indirect reference given by the First Self-Reference Theorem (and/or the Second Self-Reference Theorem) and the fixed point results of the lambda algebra available in the one categorical pair $(C,C').$
\bigbreak

\section{Self-Reference}
Finally we return to self-reference in the form of the expansion of a name. Recall the expansion operator
as described in Section 1.  We have an operation E on names that 
{\em expands}  a name in the sense that
if ``A" is the name of A then the expansion E``A" refers to A``A", the result of appending the contents of the name to the name. Thus if we regard the name as pointing to its contents as in 
$$``A" \longrightarrow A$$ then $$E``A" \longrightarrow A``A".$$ Self-reference results when one expands the name of the expansion operator. $$E``E" \longrightarrow E``E".$$ How is this self-reference
related to the self-reference we are all familiar with in our personal experience?
\bigbreak

To begin to see an answer to this question, consider the use of the pronoun ``I". When I say I then I refer to myself. I alone does not refer to itself. It is required that there be a contents related to the one who uses the word I. I am the one who says I, and this can be said by anyone.
So in a sense we can say that I is really the expansion operator and the self-reference associated with I
occurs when we apply I to ``I", forming I``I" which is self-referent. In other words, we each make a personal identification $$I = I``I",$$ that says `` I am the operation of expanding myself to
my content (which is myself)." This was said more eloquently by Heinz von Foerster
\cite{VF} : ``I am the observed relation between myself and observing myself."
We encourage the reader to expand further on these themes.
\bigbreak

\end{document}